\documentclass{article}
\usepackage[english]{babel}

\usepackage{amsmath}
\usepackage{amssymb}
\usepackage{amsfonts}
\usepackage{mathrsfs} 
\usepackage{graphics}

\newtheorem{theorem}{Theorem} 
\newtheorem{lemma}{Lemma}
\newtheorem{remark}{Remark}

\newtheorem{coro}{Corollary}

\newtheorem{defin}{Definition}

\def\{{\protect\lbrace}
\def\}{\protect\rbrace}
\def\codim{\operatorname{codim}}
\def\Real{\operatorname{Re}}
\def\rg{\operatorname{rg}}
\newcommand{\V}{\vec V}
\newcommand{\W}{\vec W}

\newcommand{\ebox}{\fbox {} \smallskip}   
\numberwithin{equation}{section}  

\title{On a class of vector fields with discontinuity of divide-by-zero type and its applications}

\author{R. Ghezzi\thanks{{\tt ghezzi@sissa.it}}, A.O. Remizov\thanks{{\tt remizov@sissa.it}}}
\date{\small{SISSA/ISAS, via Bonomea 265,
34136 Trieste,
Italy}}

\begin{document}
\large
\maketitle


\section{Introduction}

Let $\Omega$ be a domain in $\mathbb R^n$, $n \geq 2$, with coordinates $x = (x_1, \ldots,x_n)$. Let $f\in C^s(\Omega,\mathbb{R})$, $s\geq 2$, be such that the equation $f(x)=0$ defines a regular hypersurface $\Gamma \subset \Omega$, i.e., at all points $x \in \Gamma$ the condition $\nabla f(x) \neq 0$ holds.
We consider vector fields of the type
\begin{equation}
\W (x) = f^{-r}(x) \, \V (x),
\label{defW}
\end{equation}
where $\V\in C^s(\Omega,\mathbb{R}^n)$ is a vector field and $r$ is a positive real number.

The divergence of $\W$, denoted by $D_{\W}$, is infinite or undetermined on the hypersurface $\Gamma$, but it is a $C^{s-1}$-smooth function on
$\Omega \setminus \Gamma$. Assume that the field $\W$ satisfies the following main conditions
\begin{eqnarray*}
&&\lim_{x \to x_*} f^{r+1}(x) D_{\W} (x) = 0, \quad \forall\, x_* \in \Gamma,\\
&&\lim_{x \to x_*} f^r(x) D_{\W} (x) = 0, \quad \lim\limits_{x \to x_*} f^{r+1} \frac{\partial D_{\W}}{\partial x_i} (x) = 0,
\quad \forall\, x_* \in \Gamma :  \V(x_*)=0, \ \, \forall\, i.
\end{eqnarray*}
For simplicity  we  write these assumptions in the form
\begin{eqnarray}
\label{conda}
&& f^{r+1} D_{\W} {\bigl|}_{\Gamma} = 0, \\
\label{condb}
&&f^r D_{\W} (x_*) = 0, \quad f^{r+1} \frac{\partial D_{\W}}{\partial x_i} (x_*) = 0, \quad \forall\, x_* \in \Gamma :  \V(x_*)=0, \ \, \forall\, i.
\end{eqnarray}
Conditions \eqref{conda}, \eqref{condb} are fulfilled, for instance, if the vector field $\W$ is divergence-free, i.e., $D_{\W}\equiv 0$ for all $x\in \Omega \setminus \Gamma$.

The field $\W$ is  $C^s$-smooth on $\Omega \setminus \Gamma$, but at points of $\Gamma$ formula \eqref{defW} gives a discontinuity of divide-by-zero type.
Due to their large number of applications (e.g., in  mechanics with dry friction and control theory, see \cite{Fil}), discontinuous vector fields (or, equivalently, differential equations with discontinuous righthand sides) have been widely studied. However to the authors' knowledge, fields of the given type have not been studied yet.

Although  this problem seems at first sight rather theoretical and not natural, it is motivated by a large number of
applications. Indeed,  many variational problems in differential geometry and calculus of variations are
characterized by Lagrangian (or Hamiltonian) functions  that are smooth at all points except
for a regular hypersurface $\Gamma$. The vector field corresponding to the Euler--Lagrange equations of
such problems is divergence-free and takes the form \eqref{defW}. The simplest example is the equation of geodesic lines on the cuspidal edge embedded in the Euclidean space or on the plane with the Klein metric, that is used in the model of the Lobachevsky plane.

The aim of this paper is to establish some general facts about vector fields of the form \eqref{defW} under assumptions \eqref{conda}, \eqref{condb}
that allow to infer some properties on the vector fields $\V$. Applications to several concrete problems, that are interesting from different points of view, are provided.

The paper is orginized as follows. In  section~\ref{sec-2} we prove several simple theorems about vector (and direction) fields of the form \eqref{defW} under assumptions \eqref{conda}, \eqref{condb} without any special hypothesis on $\V$. In particular, these results show the key role of singular points of the field $\V$ in the applications.

In  section~\ref{sec-3} we give a brief survey of  the theory of normal forms at non-isolated singular points of smooth vector fields. We restrict to the case where the components of the vector field belong to the ideal  generated by two of them in the ring of smooth functions. As far as we know,  the first work devoted to the analysis of local normal forms for such fields is due to F.~Takens \cite{Tak}. Later, the problem was deeply investigated in finite smooth \cite{Rem06}, $C^{\infty}$-smooth \cite{Rouss}, and analytic categories \cite{Voron1,Voron2}.\footnote{
Remark that the finite smooth classification is based on the general results by V.S.~Samovol \cite{Sam}, and the analytic classification is based on the general results by A.D.~Bryuno \cite{Bryuno1}\,--\,\cite{Bryuno3} and J.C.~Yoccoz \cite{Yo}.
}

In this survey, we deal only with finite and $C^{\infty}$-smooth classifications, which are simpler than the analytic one.
Notice that almost all the facts in this section were previously known. We try to present the subject so that it is not obscured by technical details and, at the same time, is sufficiently precise. We hope that the informed reader will tolerate trivial aspects while the reader unfamiliar with
this subject will understand the main ideas and find all the omitted proofs and technical details in the
cited literature (see also \cite{IYa, Pazii, ZhitSo, Zhit} devoted to similar problems).

In the last section we apply the results to the problem of  geodesic flow generated by  three different types of singular metrics on 2-surfaces. Firstly, we consider pseudo-Riemannian metrics, i.e., metrics that degenerate (change their signature) on a curve, see also~\cite{Rem-Pseudo}.
Secondly,  we analyse  metrics of Klein type, that are positive definite but  have a singularity of divide-by-zero type, see~\cite{Rem-Klein}. Finally, we consider almost-Riemannian metrics, i.e., metrics whose orthonormal frames are pair of vector fields that are collinear on a regular curve, see~\cite{ABS}. Two another examples can be found in \cite{Rem07, Rem-Edge}.

\section{Basic Theorems}\label{sec-2}

Integral curves of the fields $\W$ and $\V$ coincide at all points $x\in \Omega \setminus \Gamma$. At the same time the field $\V$ is more suitable for analysis, since it is smooth on the whole domain $\Omega$ while the field $\W$ is discontinuous on the hypersurface $\Gamma \subset \Omega$.
Our concern is to pass from the initial vector field $\W$ to the vector field $\V$.

\begin{theorem}
\label{T1}
Condition \eqref{conda} holds true if and only if \,$\Gamma$ is an invariant hypersurface of $\V$.
The function $f$ is a first integral of $\V$ if and only if 
\begin{equation}
f^r D_{\W}(x) \equiv D_{\V}(x).
\label{firstintegral}
\end{equation}
Assume  $f$ to be a first integral of $\V$ and let condition \eqref{condb} holds true.
Then $D_{\V}(x_*) = 0$ for every $x_*\in\Gamma$ such that $\V(x_*)=0$.
\end{theorem}

{\bf Proof}.
Using the formula of  divergence in Cartesian coordinates, for every point in $\Omega\setminus\Gamma$ we get 
\begin{equation}
f^{r+1}D_{\W} = f D_{\V} + f^{r+1}L_{\V} (f^{-r}) = f D_{\V} - r L_{\V} f,
\label{diver}
\end{equation}
where $L_{\V}$ denotes the Lie derivative along the vector field $\V$. All terms in the right hand side of equation 
\eqref{diver} are $C^{s-1}$-smooth on $\Omega$. Hence, taking the limit as $x$ tends to $x_*\in\Gamma$ it follows that 
$f^{r+1} D_{\W} {\bigl|}_{\Gamma} = 0$ is equivalent to $ L_{\V} f {\bigl|}_{\Gamma} =0$.

As concerns the second statement, for every point in $\Omega\setminus \Gamma$, we have
\begin{equation}\label{diver2}
f^rD_{\W}-D_{\V}=-r f^{-1}L_{\V}f.
\end{equation}
If identity \eqref{firstintegral} holds on  $\Omega\setminus\Gamma$, then  \eqref{diver2} implies  $L_{\V}f{\bigl|}_{\Omega\setminus\Gamma} = 0$. By continuity it follows  $L_{\V}f\equiv 0$ on $\Omega$, i.e.,  $f$ is a first integral of  $\V$. Conversely, assume  $L_{\V}f \equiv 0$ on $\Omega$. Then, by \eqref{diver2} we get 
$(f^rD_{\W}-D_{\V}){\bigl|}_{\Omega\setminus\Gamma}=0$. Thus, by continuity ($D_{\V}$ is continuous on $\Omega$), identity \eqref{firstintegral} on $\Omega$ follows.
Finally, combining  the first equality in condition \eqref{condb} and identity \eqref{firstintegral}, we get the  last statement of the theorem.
\hfill\ebox

\begin{coro}
\label{C1}
Assume that condition \eqref{conda} holds. 
Let $\gamma$ be an integral curve of either $\V$ or $\W$ passing through the point $x_* \in \Gamma$.
If $\V(x_*)\neq 0$, then in a neighborhood of $x_*$ the curve $\gamma$ lies entirely in the hypersurface $\Gamma$. 
\end{coro}

Theorem \ref{T1} and  Corollary \ref{C1} explain why  singular points of  $\V$ play an important role. 
Indeed, in many applications it is necessary to find integral curves that intersect the invariant hypersurface $\Gamma$ but do not belong entirely to $\Gamma$. Hence such integral curves intersect $\Gamma$ only at singular points. The next theorem establishes a relation between the eigenvalues of the linearization of the vector field $\V$ at a singular point $x_* \in \Gamma$. As we shall see, in many cases such a relation is resonance.

%
\begin{figure}[h!]
\begin{center}
\includegraphics{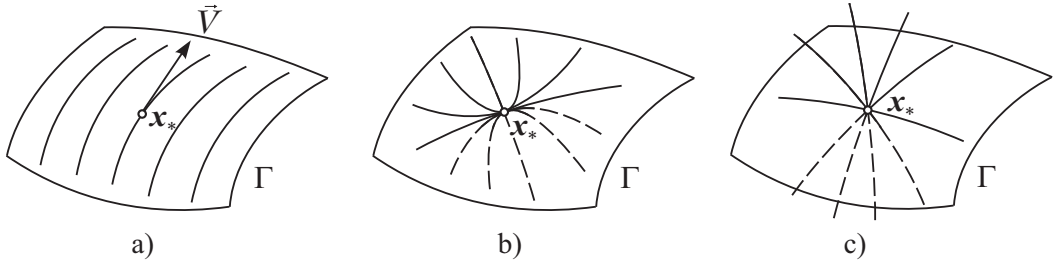}
\caption{
Three examples of phase portraits of the vector field~$\V$. a)~The~case $\V(x_*)\neq 0$, and all integral curves belong to~$\Gamma$. b),~c) The case $\V(x_*)=0$, and all integral curves (except only one) do not belong to~$\Gamma$.
}\label{fig-intro}
\end{center}
\end{figure}
\begin{theorem}
\label{T2}
Let $x_* \in \Gamma$ be a singular point of the field $\V$ and $\lambda_1, \ldots, \lambda_n$ be the eigenvalues of the linearization of $\V$ at $x_*$. 
If conditions \eqref{conda} and \eqref{condb} hold true then there exists $j \in \{1, \ldots, n\}$ such that
\begin{equation}
\lambda_1 + \cdots + \lambda_n = r\lambda_j.
\label{resonance}
\end{equation}
There exists an eigenvector corresponding to $\lambda_j$ which is transversal to $\Gamma$ at $x_*$. The spectrum of the linearization of the restriction $\V {\bigl|}_{\Gamma}$ at $x_*$ is $\{\lambda_1, \ldots, \lambda_n\} \setminus \lambda_j$.
\newline
If  $f$ is a first integral of $\V$ then $\lambda_j=0$.
\end{theorem}

{\bf Proof}. 
Differentiating  the identity $f^{r+1}D_{\W}=fD_{\V}-rL_{\V}f$ (see \eqref{diver}) with respect to  $x_i$ we get
\begin{equation*}
(r+1)f^r D_{\W} \, \frac{\partial f}{\partial x_i} + f^{r+1} \frac{\partial D_{\W}}{\partial x_i} = D_{\V} \, \frac{\partial f}{\partial x_i} + f \frac{\partial D_{\V}}{\partial x_i} - r \biggl\langle \frac{\partial \V}{\partial x_i}, \nabla f\biggr\rangle - r \biggl\langle \V, \nabla \frac{\partial f}{\partial x_i} \biggr\rangle,
\end{equation*}
where the triangle brackets denote the standard scalar product of vectors.

The last equality holds for all $x \in \Omega \setminus \Gamma$ and its right hand side is $C^{s-2}$-smooth on $\Omega$.
Taking  the limit as $x$ tends to $x_* \in \Gamma$ such that $\V(x_*)=0$ and using \eqref{condb}, we get
$$
\left(D_{\V} \, \frac{\partial f}{\partial x_i} - r \biggl\langle \frac{\partial \V}{\partial x_i}, \nabla f\biggr\rangle\right){\biggl|}_{x_*} = 0, \quad i=1, \ldots, n.
$$
Since $r \neq 0$, the last system  can be written in  matrix form as 
$A e = \rho e$, where $A = \bigl( \frac{\partial \V}{\partial x}\bigr) {\bigl|}_{x_*}$ is the matrix of the linearization of $\V$ at the singular point $x_*$, the vector $e = \nabla f(x_*)$, and the number $\rho = r^{-1} D_{\V}(x_*)$. 
By hypothesis,  $\nabla f(x_* ) \neq 0$, hence $\rho$ is an eigenvalue of the linearization of $\V$ at $x_*$ with corresponding eigenvector $\nabla f(x_*)$. Let $j\in\{1,\dots n\}$ be such that $\rho=\lambda_j$.
Notice that the divergence of a vector field at any singular point coincides with the trace of the linearization of this field at that point. Thus $\lambda_j=r^{-1}(\lambda_1 + \cdots + \lambda_n)$ which leads  to equality \eqref{resonance}. Clearly, $e$ is transversal to $\Gamma$ at $x_*$, whence  $\lambda_j$ does not belong to the spectrum of the linearization of the restriction $\V {\bigl|}_{\Gamma}$. 

To prove the last statement  recall that if  $f$ is a first integral of $\V$, by Theorem \ref{T1}  
$D_{\V}(x_*) = 0$  follows. Since $D_{\V}(x_*) = \lambda_1 + \cdots + \lambda_n$, we have equality \eqref{resonance} with $\lambda_j=0$.
\hfill\ebox

\begin{remark}
\label{R1}
Theorems \ref{T1} and \ref{T2} hold true not only for $r>0$, but also for $r<0$. 
\end{remark}

Let us illustrate two examples in $\mathbb{R}^3$ with coordinates $(x,y,z)$.

{\bf Example 1}. Let $\V$ be the vector field
\begin{equation*}
\dot x= x,\quad \dot y= y, \quad \dot z=z.
\end{equation*}
Then $D_{\V}(x,y,z) \equiv 3$ and the unique singular point of $\V$ is the origin. The spectrum of the linearization at the origin is $(1,1,1)$. The field $\V$ has no non-constant first integrals, but it has the family of integral planes $ax+by+cz=0$ passing through the origin.

Consider the vector field $\W$ given by formula \eqref{defW} with $f(x,y,z) = ax+by+cz$. Then $D_{\W} = (3-r) f^{-r}$ and $f^r D_{\W} = (3-r)$. Condition \eqref{conda} is satisfied, but \eqref{condb} is fulfilled only if $r=3$. In the case $r=3$, we have relation \eqref{resonance} with any index $j=1,2,3$.

{\bf Example 2}. Let $\V$ be the vector field
\begin{equation*}
\dot x= 2x,\quad \dot y= y, \quad \dot z=0.
\end{equation*}
Then  $D_{\V}(x,y,z) \equiv 3$ and the set of singular points of $\V$ is the $z$-axis. The spectrum of the linearization at any singular point 
is $(2,1,0)$. The coordinate function $z$ is a first integral of $\V$ and there is a family of integral surfaces given by  $x - cy^2 = 0$, as $c$ varies in $\mathbb{R}$. 

Consider the vector field $\W$ given by formula \eqref{defW} with $f(x,y,z)=z$. Then $f^r D_{\W} = 3$ and condition \eqref{condb} is not satisfied.
This corresponds to the last claim of Theorem~\ref{T1}. Indeed, since  $f$ is a first integral of $\V$, condition \eqref{condb} would imply $D_{\V}(0,0,z)\equiv 0$.

Let now  $f(x,y,z) = x - cy^2$.  Defining $\W$ as in \eqref{defW}, we get $f^r D_{\W} = (3-2r)$. Hence condition \eqref{condb} holds true if and only if $r=3/2$ and the  relation  \eqref{resonance} is valid with  $\lambda_j=2$. 

Finally, consider the vector field $\W$ with $f(x,y,z) = y$. Then $f^r D_{\W} = (3-r)$, condition \eqref{condb} holds in the case $r=3$, and we have \eqref{resonance} with $\lambda_j=1$.

\medskip 

Sometimes it is more natural to consider direction fields rather than vector fields. Recall that given a vector field $\V$, the  direction field $\chi$ associated to $\V$ is the class of vector fields $\varphi \V$, where $\varphi\in C^s(\Omega)$ never vanishes. Theorems \ref{T1}, \ref{T2} are valid for direction fields.

\begin{theorem}
\label{T3}
Let $\varphi\in C^{s}(\Omega)$ and $\varphi(x)\neq 0$ for every $x\in\Omega$. 
Then Theorems \ref{T1}, \ref{T2} hold true if in equation \eqref{defW} we replace $\V$ by $\varphi\V$.
\end{theorem}

{\bf Proof}.
It is necessary and sufficient to prove that the main assumptions \eqref{conda} and \eqref{condb} are invariant with respect to multiplication of the vector fields $\V$ (and consequently $\W$) by a $C^s$-smooth scalar factor $\varphi \neq 0$. Indeed, $D_{\varphi \W} = \varphi D_{\W} + f^{-r} L_{\V} \varphi$.
Hence we get
\begin{align*}
& f^r D_{\varphi \W}{\Bigl|}_{\Gamma}= (\varphi f^{r}D_{\W}  + L_{\V} \varphi){\Bigl|}_{\Gamma}, \\
& f^{r+1} \frac{\partial D_{\varphi \W}}{\partial x_i} {\Bigl|}_{\Gamma} =\left(
f^{r+1} \Bigl( \frac{\partial \varphi}{\partial x_i} D_{\W} + \varphi \frac{\partial D_{\W}}{\partial x_i} \Bigr) + f \frac{\partial L_{\V}\varphi}{\partial x_i} 
- r \frac{\partial f}{\partial x_i} L_{\V}\varphi \right) {\Bigl|}_{\Gamma}.
\end{align*}
These expressions show that conditions \eqref{conda}, \eqref{condb} hold true for the vector fields $\varphi\V$, $\varphi\W$.

\hfill\ebox

\section{Fields with non-isolated singular points}\label{sec-3}

We are interested in studying vector fields $\V$ of the form
\begin{equation}
\dot \xi = v, \quad \dot \eta = w, \quad \dot \zeta_i = \alpha_i v + \beta_i w, \quad i=1,\ldots, l,
\label{16}
\end{equation}
where $\alpha_i, \beta_i$ and $v, w$ are $C^{\infty}$-smooth functions of the variables $\xi,\eta, \zeta_1,\ldots,\zeta_l$.
Such a kind of fields occurs in many problems, for instance, 
in studying implicit differential equations (see next example) and slow-fast systems.

{\bf Example 3}.
Consider the family of first-order implicit differential equations 
\begin{equation}
\label{esem4}
F(t,x,p)=\varepsilon,  \quad p=\frac{dx}{dt},
\end{equation}
depending on the real parameter $\varepsilon$ not necessarily small.
One effective approach (which goes back to Poincar\'e) consists of lifting the multi-valued direction field defined by equation \eqref{esem4} on the $(t,x)$-plane to a single-valued direction field $\chi$ defined by equation \eqref{esem4} in the $(t,x,p)$-space.

Geometrically, $\chi$ is the intersection between the contact planes $dx = pdt$ and the planes tangent to the
 surfaces $\{F=\varepsilon\}$ with various $\varepsilon$. This gives the Pfaffian system
\begin{equation*}
F_t \, dt + F_x \, dx + F_p \, dp = 0, \quad  p \, dt - dx = 0.
\end{equation*}
Whence the direction field $\chi$ corresponds to the vector field $\V$ given by the formula
\begin{equation}
\dot t = F_p, \quad \dot x = pF_p, \quad \dot p = -(F_t + pF_x),
\label{implode}
\end{equation}
where a dot over a symbol means differentiation with respect to the independent variable playing the role of time. 
The field \eqref{implode} has the form \eqref{16} with $l=1$, where $\xi=t$, $\eta=p$, $\zeta=x$ and $v=F_p$, $w=-(F_t + pF_x)$.
\medskip

In this section we recall smooth local normal forms of fields \eqref{16} at singular points.

The components of  $\V$  belong to the ideal $I=(v,w)$ generated by two of them in the ring of  germs of $C^{\infty}$-smooth functions (this property is invariant with respect to the action of diffeomorphisms of the phase space). The set of singular points of  $\V$ is defined by the two equations $v=w=0$. The spectrum of the linearization of $\V$ at any singular point contains at least $l$ zero eigenvalues, i.e., it is $(\lambda_1, \lambda_2, 0, \ldots, 0)$.

Consider the germ of  $\V$ at a given  singular point. Without loss of generality we may assume the singular point to be  the origin of the phase space. From now on,  we will always assume that  $\Real \lambda_{1,2}(0) \neq 0$, whence the set of singular points of $\V$ is the regular center manifold of dimension $l$,  denoted by $W^c$. The eigenvectors with zero eigenvalue are tangent to $W^c$ and the eigenvectors (if they exist) corresponding to $\lambda_{1,2}(0)$  are tangent to the plane $d\zeta_i = \alpha_i \, d\xi + \beta_i \, d\eta$, $i=1,\ldots, l$.

It is convenient to choose  local coordinates $(\xi,\eta, \zeta_1,\ldots,\zeta_l)$ such that $W^c=\{\xi=\eta=0\}$ and the linear
part  of $\V$ at $0$ is in normal Jordan form. Then there exist $C^{\infty}$-smooth functions $v_{1,2}$ and $w_{1,2}$ such that $v = \xi v_1 + \eta v_2$ and  $w = \xi w_1 + \eta w_2$. The eigenvalues $\lambda_{1,2}$ at various singular points continuously depend on the variable $\zeta = (\zeta_1,\ldots,\zeta_l)$, which is a local coordinate on $W^c$ (this dependence is $C^{\infty}$-smooth at the points where $\lambda_1 \neq \lambda_2$). In the following we will always work using such coordinates.

\subsection{Normal forms: the non-resonant case}

We shall say that $k$ functions $U^{(i)}(\xi,\eta,\zeta)$, $i=1,\ldots,k$, are {\it independent by} $\zeta$ at the point $0$ if their gradients with respect to the variable $\zeta = (\zeta_1,\ldots,\zeta_l)$ at $0$ are linearly independent. It is not hard to see that if a function $U$ is a first integral of  $\V$ then its partial derivatives  $U_{\xi}$ and $U_{\eta}$ vanish at $0$. Hence the number of first integrals of $\V$ independent at $0$ is not greater than $l$. On the other hand, the restriction of $\V$ to  $W^c$ is identically  zero, hence by Shoshitaishvili's reduction theorem \cite{AI},
the germ of $\V$ at $0$ is orbitally topologically equivalent to
\begin{equation*}
\dot \xi =  \xi, \quad \dot \eta = \pm \eta, \quad \dot \zeta_i = 0,  \quad i=1,\ldots, l.
\end{equation*}

The trivial equations $\dot \zeta_i = 0$ suggest the existence of $l$ independent first integrals  of $\V$.
If a $l$-uple of smooth first integrals $U^{(1)}, \ldots, U^{(l)}$ independent by $\zeta$ at $0$ exists,
the change of coordinates $\zeta_i \mapsto U^{(i)}$ brings the field $\V$ to the form
\begin{equation*}
\dot \xi = \overline v, \quad \dot \eta = \overline w, \quad \dot \zeta_i = 0,  \quad i=1,\ldots, l,
\end{equation*}
where $\overline v$ and $\overline w$ are smooth functions obtained from $v$ and $w$ by the above change of coordinates.

The existence of $l$ independent smooth first integrals is connected with the following concept.

\begin{defin}
\label{D1}
The relations
\begin{equation}
p_1\lambda_1 + p_2\lambda_2 = 0, \quad p_{1,2} \in \mathbb Z_+, \quad p_1+p_2 \geq 1,
\label{18}
\end{equation}
are called resonances of first type.
The minimal number $p_1+p_2$ (i.e., $p_1$ and $p_2$ are relatively prime) is called the order of the resonance \eqref{18}.
\end{defin}

Consider the germ of a smooth function $U(\xi,\eta,\zeta)$ at the point $0$ and its Taylor series with respect to the variables $\xi,\eta$, i.e.,
\begin{equation}
U(\xi,\eta,\zeta) = \sum_{p_{1,2} \in \mathbb Z_+}^{} u_{p_1p_2}(\zeta)\, \xi^{p_1} \eta^{p_2}.
\label{19}
\end{equation}
The germ of $U$ is called {\it $N$-flat} ($N \in \mathbb N$ or $\infty$) by $\xi,\eta$ if $u_{p_1p_2}(\zeta) \equiv 0$ for all $p_1+p_2 \leq N$.

\begin{lemma}
\label{L1}
If between the eigenvalues $\lambda_{1,2}(0)$  there are no resonances \eqref{18} up to  order $N \in \mathbb N$ inclusive,
then there exist $C^{\infty}$-smooth functions $U^{(1)}, \ldots, U^{(l)}$ independent by $\zeta$ at $0$
such that  $L_{\V}U^{(1)}, \ldots, L_{\V}U^{(l)}$ are $N$-flat by $\xi,\eta$ at $0$.
\end{lemma}

{\bf Proof}.
We prove the lemma assuming $\lambda_{1,2}(0)$ to be real and the Jordan form of the linearization of $\V$ at $0$ to be diagonal (the cases of complex eigenvalues or Jordan form with a second-order cell to be considered similarly). Then the germ of $\V$ has the form
\begin{equation}
\label{20}
\begin{aligned}
\dot \xi &= \xi(\lambda_1(0)+{\widetilde v}_1(\zeta)+\cdots) + \eta({\widetilde v}_2(\zeta)+\cdots), \\
\dot \eta &= \xi({\widetilde w}_1(\zeta)+\cdots) + \eta(\lambda_2(0)+{\widetilde w}_2(\zeta)+ \cdots), \\
\dot \zeta_i &= \xi(h_1^{(i)}(\zeta)+\cdots) + \eta(h_2^{(i)}(\zeta)+\cdots),   \quad \   i=1,\ldots, l, \\
\end{aligned}
\end{equation}
where all functions ${\widetilde v}_{1,2}(\zeta)$, ${\widetilde w}_{1,2}(\zeta)$, $h_{1,2}^{(i)}(\zeta)$ are $C^{\infty}$-smooth and vanish at $0$, and the omitted terms are $C^\infty$-smooth functions containing the factor $\xi$ or $\eta$.

The idea is to look for the functions $U^{(1)}, \ldots, U^{(l)}$ in the set of polynomials in $\xi,\eta$ with coefficients smoothly depending on $\zeta$. 
Namely, consider a function $U$ in the form \eqref{19} with finite sum $0 \leq p_1+p_2 \leq N$ and unknown coefficients $u_{p_1 p_2}(\zeta)$. Substituting this expression into $L_{\V}U$ and using \eqref{20}, we get
\begin{multline*}
L_{\V}U = \sum_{p_1+p_2=0}^{N} \biggl(
u_{p_1 p_2} \xi^{p_1}\eta^{p_2} \bigl(p_1(\lambda_1(0)+{\widetilde v}_1+\cdots)+p_2(\lambda_2(0)+{\widetilde w}_2+\cdots)\bigr) + \\
+ u_{p_1 p_2} p_1 \xi^{p_1-1}\eta^{p_2+1} ({\widetilde v}_2+\cdots) + u_{p_1 p_2} p_2 \xi^{p_1+1}\eta^{p_2-1} ({\widetilde w}_1+\cdots) + \\
+ \sum_{k=1}^l  \frac{\partial u_{p_1p_2}}{\partial \zeta_k} \Bigl(\xi^{p_1+1}\eta^{p_2}(h_1^{(k)}+\cdots) + \xi^{p_1}\eta^{p_2+1}(h_2^{(k)}+\cdots)\Bigr) \biggr).
\end{multline*}

For $L_{\V} U$ to be $1$-flat by $\xi, \eta$ we set the coefficients of the monomials $\xi$ and $\eta$ equal to zero, that is,
\begin{equation}
\label{21}
\begin{aligned}
(\lambda_1(0)+{\widetilde v}_1) u_{10} + {\widetilde w}_1 u_{01} + \sum_{k=1}^l h_1^{(k)} \frac{\partial u_{00}}{\partial \zeta_k} = 0, \\ 
{\widetilde v}_2 u_{10} + (\lambda_2(0)+{\widetilde w}_2) u_{01} + \sum_{k=1}^l h_2^{(k)} \frac{\partial u_{00}}{\partial \zeta_k} = 0. \\
\end{aligned}
\end{equation}
Since the determinant $d_1(\zeta)$ of the linear system \eqref{21} with respect to the unknown variables $u_{10}$ and $u_{01}$ is a smooth function and $d_1(0) = \lambda_1(0) \lambda_2(0) \neq 0$, the solutions $u_{10}$ and $u_{01}$ are smooth in a neighborhood of $0$. Notice that the functions $u_{10}(\zeta)$ and $u_{01}(\zeta)$ depend on the derivatives $\frac{\partial u_{00}}{\partial \zeta_i}$, where $u_{00}(\zeta)$ is any arbitrary smooth function.

Given $n \in \{2, \ldots, N\}$ consider the coefficients of the monomials $\xi^{p_1}\eta^{p_2}$, where $p_1+p_2=n$.
In order $L_{\V} U$ to be $n$-flat, we set these coefficients to be identically zero, i.e., we solve the system 
\begin{equation}
(p_1(\lambda_1(0)+{\widetilde v}_1)+p_2(\lambda_2(0)+{\widetilde w}_2)) \, u_{p_1p_2} + 
(p_1+1){\widetilde v}_2 \, u_{p_1+1 p_2-1} + (p_2+1){\widetilde w}_1 \, u_{p_1-1 p_2+1} = \varphi_{p_1p_2},
\label{22}
\end{equation}
where $\varphi_{p_1p_2}$ are polynomials of the coefficients $u_{\alpha \beta}(\zeta)$ and their first-order derivatives with $\alpha+\beta < n$ 
(with $u_{\alpha \beta}(\zeta) \equiv 0$ if $\alpha<0$ or $\beta<0$).
The determinant $d_n(\zeta)$ of the linear system \eqref{22} with respect to  variables $u_{p_1p_2}$, $p_1+p_2=n$, has the form
\begin{equation*}
d_n(\zeta) = \prod\limits_{p_1+p_2=n} (p_1\lambda_1(0) + p_2\lambda_2(0) + \delta(\zeta)),
\end{equation*}
where $\delta(\zeta)$ is a smooth function vanishing at $\zeta=0$. The absence of resonances \eqref{18} up to order $N$ implies that $d_n(0) \neq 0$, whence in a neighborhood of $0$ the coefficients $u_{p_1p_2}(\zeta)$, $p_1+p_2=n$, smoothly depend on the functions $u_{\alpha \beta}(\zeta)$ and their first-order derivatives with $\alpha+\beta < n$.

Finally, let $u_{00}^{(1)}(\zeta),\dots, u_{00}^{(l)}(\zeta)$ be $C^{\infty}$-smooth functions independent at $0$. For every index $i=1,\dots, l$ we define $U^{(i)}$ by solving systems \eqref{21} and \eqref{22} for $n = 2, \ldots, N$ with the initial function $u_{00}=u_{00}^{(i)}$. By construction, $U^{(1)},\dots, U^{(l)}$ satisfy the required conditions.
\hfill\ebox

\begin{coro}
\label{C2}
If between the eigenvalues $\lambda_{1,2}(0)$ there are no resonances \eqref{18} up to  order $N$ inclusive, then the last $l$ components 
$\alpha_i v + \beta_i w$ of the field \eqref{16} can be assumed to be $N$-flat by $\xi,\eta$ at $0$.
\end{coro}

Corollary \ref{C2} allows to get normal forms in $C^{k}$-smooth and $C^{\infty}$-smooth categories.
As for the $C^k$-smooth category, in \cite{Sam} the author defines the number
\begin{equation*}
N(k) = 2 \biggl[ (2k+1) \frac{\max |\Real \lambda_{1,2}|}{\min \, |\Real \lambda_{1,2}|} \, \biggr] + 2,  \quad k \in \mathbb N,
\end{equation*}
the square brackets denoting the integer part of a number.

\begin{theorem}
\label{Ck}
If between the eigenvalues $\lambda_{1,2}(0)$ there are no resonances \eqref{18} up to  order $N(k)$, 
then the germ of \eqref{16} is $C^k$-smoothly equivalent to 
\begin{equation}
\label{cknormalform}
\dot \xi = v, \quad \dot \eta = w, \quad \dot \zeta_i=0, \quad i=1,\ldots, l,
\end{equation}
where $v,w$ are some new functions of $\xi,\eta,\zeta$.
If between the eigenvalues $\lambda_{1,2}(\zeta)$ there are no resonances \eqref{18} of any order for all $\zeta$ sufficiently close to $0$, 
then the germ of \eqref{16} is $C^\infty$-smoothly equivalent to
\eqref{cknormalform}.
\end{theorem}

The proof of Theorem \ref{Ck} in the finite-smooth category is based on Lemma \ref{L1} and on general results from \cite{Sam}.
The proof in the $C^\infty$-smooth category requires more advanced techniques (see \cite{Rouss} or \cite{IYa}).
Notice that if $\Real \lambda_{1,2}(0)$ have the same sign the absence of resonances \eqref{18} between $\lambda_{1,2}(0)$ implies the absence of resonances \eqref{18} between $\lambda_{1,2}(\zeta)$ for all $\zeta$ sufficiently close to $0$. This is no longer true if $\Real \lambda_{1,2}(0)$ have different signs, except for the special case when the ratio $\lambda = \lambda_{1}/\lambda_2$ is constant on $W^c$, i.e., at all singular points.

As we shall see in the following, the normal form \eqref{cknormalform} can be further simplified. 
\begin{defin}
The relations
\begin{equation}
p_1 \lambda_1 + p_2 \lambda_2 = \lambda_j, \quad p_{1,2} \in \mathbb Z_+, \quad p_1+p_2 \geq 2, \quad j \in \{1,2\},
\label{24}
\end{equation}
are called resonances of second type. The number $p_1+p_2$ is the {\it order} of  resonance. 
\end{defin}

Clearly, a resonance \eqref{18} of order $n$ implies a resonance \eqref{24} of order $n+1$. In this section we assume the absence of resonances \eqref{18} up to order $N \in \mathbb N$ or $\infty$. Hence  a resonance \eqref{24} of order $\leq N$ holds if and only if the ratio $\lambda(0)=\lambda_1(0)/\lambda_2(0)$ or its inverse belongs to $\{2,\ldots, N\}$.
Combining the results from \cite{IYa}, \cite{Rem06}, \cite{Rouss}, \cite{Sam}, one gets the following theorems.

\begin{theorem}
\label{T2-1}
Let $k\in\mathbb{N}$ and assume that between $\lambda_{1,2}(0)$ there are no resonances \eqref{24} of order $N(k)$ inclusive.
Then the germ of $\V$ at $0$ is $C^k$-smoothly equivalent to
\begin{equation}
\dot \xi = \alpha_1(\zeta)\xi + \alpha_2(\zeta)\eta, \quad \dot \eta = \alpha_3(\zeta)\xi + \alpha_4(\zeta)\eta, \quad \dot \zeta_i = 0, \quad i=1,\ldots, l.
\label{25}
\end{equation}
Moreover, if $\lambda_{1,2}(0)$ are real and $\lambda(0)\neq 1$,  the germ of $\V$ at $0$ is $C^k$-smoothly orbitally equivalent to
\begin{equation}
\dot \xi = \lambda(\zeta)\xi, \quad \dot \eta = \eta, \quad \dot \zeta_i = 0,  \quad i=1,\ldots, l.
\label{26}
\end{equation}
Both statements hold true with $k=\infty$ if between $\lambda_{1,2}(\zeta)$ there are no resonances \eqref{24} of any order for all $\zeta$  sufficiently close to $0$.
\end{theorem}

\begin{theorem}
\label{T2-3}
Assume that $\lambda(0)=n$ is natural. Then the germ of $\V$ at $0$ is $C^{\infty}$-smoothly orbitally equivalent to 
\begin{equation}
\dot \xi = \lambda(\zeta)\xi + \varphi(\zeta) \eta^n, \quad \dot \eta = \eta, \quad \dot \zeta_i = 0,  \quad i=1,\ldots, l.
\label{27}
\end{equation}
If $\varphi(0) \neq 0$, then $\varphi(\zeta)$ simplifies to~$1$; if $\varphi(\zeta)$ has a zero of finite order $s$ at the origin then $\varphi(\zeta)$ simplifies to $\zeta^{s}$.
\end{theorem}

The normal forms \eqref{25}--\eqref{27} show that in a small neighborhood of $0$ the phase portrait of $\V$ is rather simple and $\V$ has a smooth 2-dimensional invariant foliation given by the equation $\zeta=c$ in normal coordinates.
The restriction of $\V$ to each leaf $\zeta=c$ is a planar vector field with non-degenerate singular point: node, saddle, or focus.

\subsection{Normal forms: the resonant case}

Consider the case where between $\lambda_{1,2}(0)$ there is a resonance \eqref{18}, i.e., there exist $n,m\in{\mathbb N}$ relatively prime such that
\begin{equation}
m\lambda_1 + n\lambda_2 = 0.
\label{28}
\end{equation}
In this case, the proof of Lemma~\ref{L1} for $N \geq n+m$ fails, since the determinant $d_{n+m}(\zeta)$ of the linear system \eqref{22} with $p_1+p_2=n+m$ vanishes at $\zeta=0$, and vector field \eqref{16} with resonance \eqref{28} at $0$ may not have a $l$-uple of smooth first integrals independent at $0$.

A simple illustration (with $l=1$ and $n=m=1$) comes from Example~3. Indeed, let $0$ be a singular point of the vector field $\V$ given by \eqref{implode}. Clearly, $F$ is a first integral of $\V$, and the derivatives $F_p$ and $F_t$ vanish at $0$. Assume that $\lambda_{1}(0)+\lambda_{2}(0)=0$. Since $\lambda_1+\lambda_2=D_{\V}=-F_x$, we have $F_x(0)=0$, i.e., $F$ is not regular at $0$. Let $\widehat F$ be another first integral of $\V$.
Then the integral curves of $\V$ are 1-graphs of solutions of the implicit equation $\widehat F(x,y,p) = \varepsilon$ with various $\varepsilon$. Hence, the previous argument with $\widehat F$ replacing $F$ leads to the same conclusion. Thus the germ of $\V$ at $0$ admits no regular first integrals.

\medskip
The resonance \eqref{28} generates two infinite sequences of resonances \eqref{24}, namely,
$$
(mj+1)\lambda_1 + nj\lambda_2 = \lambda_1, \quad mj\lambda_1 + (nj+1)\lambda_2 = \lambda_2, \quad j\in\mathbb{N}. 
$$
This suggests that the formal normal form contains infinite number of terms $\rho^j$, $\xi \rho^j$, $\eta \rho^j$, where $\rho = \xi^m \eta^n$
is called {\it resonance monomial} corresponding to \eqref{28}.
The central step in the derivation of normal forms in the resonant case is  the following.

\begin{lemma}
\label{L2}
 For any $k \in \mathbb N$, the germ of $\V$ at $0$ is $C^k$-smoothly equivalent to
\begin{equation}
\dot \xi = \xi(\lambda_1(0)+\Phi_1(\rho,\zeta)),  \quad \dot \eta = \eta(\lambda_2(0)+\Phi_2(\rho,\zeta))_, \quad
\dot \zeta_i = \rho \Psi_i(\rho,\zeta),    \quad i=1,\ldots, l,
\label{29}
\end{equation}
where $\Phi_{1,2}(\rho,\zeta)$ and $\Psi_i(\rho,\zeta)$ are polynomials in $\rho = \xi^m \eta^n$ of degrees $N(k)$ and $N(k)-1$, respectively, 
with coefficients smoothly depending on $\zeta$.

Assume that $\Psi_1(0,0) \neq 0$. Then for any $\omega_1, \ldots, \omega_l\in \mathbb{R}$ the germ \eqref{29} has a smooth first integral $U(\rho,\zeta)$ such that
\begin{eqnarray}
&&\Phi(\rho,\zeta) U_{\rho} + \Psi_1(\rho,\zeta) U_{\zeta_1} + \cdots + \Psi_l(\rho,\zeta) U_{\zeta_l} = 0 
\label{30}\\
&&U_{\rho}(0,0)=\omega_1, \quad U_{\zeta_2}(0,0)=\omega_2, \, \ldots, \, U_{\zeta_l}(0,0)=\omega_l, 
\label{31}
\end{eqnarray}
where $\Phi(\rho,\zeta) = m\Phi_1(\rho,\zeta) + n\Phi_2(\rho,\zeta)$.
\end{lemma}

The proof of Lemma \ref{L2} is based on the general results in~\cite{Sam} and can be found in \cite{Rem06}.

From now on, we will always assume  $\Psi_1(0,0) \neq 0$.
This hypothesis implies the existence of  $l-1$ independent first integrals $U^{(2)}, \ldots, U^{(l)}$ 
given by solutions of \eqref{30} with initial conditions \eqref{31} corresponding to linearly independent $(l-1)$-uples $(\omega_2, \ldots, \omega_l)$.
Applying the change of coordinates $\zeta_i \mapsto U^{(i)}$, $i=2,\ldots, l$, the vector field \eqref{29} is $C^\infty$-smoothly equivalent to 
\begin{equation}
\dot \xi = \xi(\lambda_1(0)+\Phi_1(\rho,\zeta)),  \quad \!
\dot \eta = \eta(\lambda_2(0)+\Phi_2(\rho,\zeta))_, \quad \!
\dot \zeta_1 = \rho \Psi_1(\rho,\zeta), \quad \! \dot \zeta_i = 0, \quad \! i=2,\ldots, l,
\label{32}
\end{equation}
where $\Phi_{1,2}(\rho,\zeta)$ and $\Psi_1(\rho,\zeta)$ are smooth functions of $\rho$ and $\zeta$ 
(not necessarily polynomials in $\rho$ like in \eqref{29}), $\Phi_{1,2}(0,0)=0$, and $\Psi_1(0,0)\neq 0$.

The first integral $U(\rho,\zeta)$ given by the solution of \eqref{30} with  initial conditions $\omega_1=1, \omega_2=\dots=\omega_l=0$ allows to simplify the form \eqref{32}. Considering the restriction $\Phi(\rho,\zeta)|_{W^c}=\Phi(0,\zeta)$, we analyse two cases: $\Phi_{\zeta_1}(0,0) \neq 0$, which is generic, or $\Phi(0,\zeta) \equiv 0$, which occurs in the analysis of some concrete problems (for instance, when $n=m=1$, this condition corresponds to
divergence-free fields).

In the first case, there exists a $C^\infty$-smooth  change of coordinates that preserves the form \eqref{32} and brings the first integral satisfying \eqref{30} with initial conditions $\omega_1=1, \omega_2=\dots=\omega_l=0$ to the form $U(\rho,\zeta)=\rho + \zeta^2$. Even if the form \eqref{32} cannot be further simplified, the phase portrait of $\V$ can be described using the invariant foliation $\rho + \zeta^2 = c$, see \cite{Rem06}.

Similarly, in the second case there exists a $C^\infty$-smooth change of coordinates that preserves the form \eqref{32} and brings the first integral  satisfying \eqref{30} with initial conditions $\omega_1=1, \omega_2=\dots=\omega_l=0$ to the form $U(\rho,\zeta)=\rho$. Using this fact, the normal form \eqref{32} simplifies as follows.

\begin{theorem}
\label{T2-4}
If  conditions $\Psi_1(0,0) \neq 0$ and $\Phi(0,\zeta) \equiv 0$ in \eqref{29} hold, then the germ of $\V$ at $0$ is $C^\infty$-smoothly orbitally equivalent to
\begin{equation}
\dot \xi  = n\xi, \quad  \dot \eta = -m\eta, \quad \dot \zeta_1 = \rho, \quad \dot \zeta_i = 0, \quad \! i=2,\ldots, l.
\label{33}
\end{equation}
\end{theorem}

The normal form \eqref{33} with any $n,m \in \mathbb N$ was established in the $C^k$-smooth category for arbitrary $k \in \mathbb N$ in \cite{Rem06}. 
It was previously proved by R.~Roussarie for the partial case $n=m=1$ in $C^\infty$-smooth category \cite{Rouss}. The techniques developed in \cite{Rouss} can be applied to establish the normal form \eqref{33} with any $n,m \in \mathbb N$ in the $C^\infty$-smooth category. However, to the authors' knowledge, this result is not published.

\begin{remark}
\label{R2}
Theorem \ref{T2-4} is not valid in the analytic case: the analytic normal form is obtained from the smooth normal form \eqref{33} by adding some module, see \cite{Voron1,Voron2}. From the viewpoint of the general theory developed by A.\,D.~Bryuno \cite{Bryuno1}\,--\,\cite{Bryuno3}, this can be explained in the following way. The \textit{condition $A$} for the formal normal form \eqref{33} does not hold, since the third equation in \eqref{33}
has the form $\dot \zeta_1 = \rho$ instead of $\dot \zeta_1 = 0$. Moreover, the pair of nonzero eigenvalues $\lambda_{1,2}(0)$ lies in the Siegel domain, where formal normalizing series generally diverge.
\end{remark}

\begin{remark}
\label{R3}
The condition $\Psi_1(0,0) \neq 0$ in Lemma \ref{L2} can be replaced by $\Psi_i(0,0) \neq 0$ for some  $i \in \{1,\ldots,l\}$. 
This condition holds true for germs \eqref{16} with resonance \eqref{28} having generic $(n+m)$-jet. Moreover, in order to check this condition it is sufficient to bring only the $(n+m)$-jet of \eqref{16} to the form \eqref{29}.
\end{remark}

The following example shows that the condition $\Psi_i(0,0) \neq 0$ is essential. 

{\bf Example 4.}
Consider the vector fields
\begin{eqnarray}
&& \dot \xi  = \xi,  \quad \dot \eta = -\eta,  \quad \dot \zeta = 0,
\label{34000}\\
&& \dot \xi = \xi,  \quad \dot \eta = -\eta (1+\xi\eta),  \quad \dot \zeta = \xi\eta\zeta,
\label{35000}
\end{eqnarray}
both having at each singular point the resonance \eqref{28} with $n=m=1$, whence  $\Phi(0,\zeta) \equiv 0$.
Clearly, the plane $\{\zeta=0\}$ is invariant for both the vector fields and it is transversal to the center manifold $W^c=\{\xi=\eta=0\}$ at the origin. If the germ of either \eqref{34000} or \eqref{35000} were $C^k$-smoothly ($k \geq 2$) orbitally equivalent to the normal form \eqref{33}, then  \eqref{33} had a $C^k$-smooth invariant surface transversal to the $\zeta$-axis, i.e., of the form $\zeta = f(\xi,\eta)$. On the other hand,  substituting the Taylor expression (of the second degree) of the function $f(\xi,\eta)$ into equation $\xi f_{\xi} - \eta f_{\eta} - \xi\eta = 0$, it is not hard to see that \eqref{33} can not have an invariant surface of the form $\zeta = f(\xi,\eta)$.

\begin{remark}
\label{R4}
If $n+m$ is rather large and the ratio $n/m$ is sufficiently close to $1$, the inequality $n+m > 2[(2k+1)\max\,\{n/m, m/n\}]+2$ has solutions 
$k \in \mathbb{N}$. According to Theorem~\ref{T2-1}, for any such $k$ the germ of \eqref{33} is $C^k$-smoothly orbitally equivalent to \eqref{26} with $\lambda(\zeta)\equiv -n/m$ or, equivalently, to the field
$$
\dot \xi  = n\xi, \quad \dot \eta = -m\eta, \quad \dot \zeta_i = 0, \quad  i=1,\ldots, l.
$$
\end{remark}

\section{Applications: geodesic flows on surfaces with singular metrics}

We start with some general consequences of the results in two previous sections and then apply them to several concrete problems connected with singularities of divide-by-zero type.

Let $\W$ be a vector field of the type in \eqref{defW}, where $r \neq 0, 1$ and the smooth\footnote{
For simplicity, we always assume that smooth means $C^{\infty}$-smooth.
}
vector field $\V$ has the form \eqref{16}. Assume that conditions \eqref{conda} and \eqref{condb} hold true.
Let $0$ be a singular point of $\V$ such that the linearization of $\V$ at $0$ has at least one non-zero real eigenvalue, i.e.,
the spectrum is $(\lambda_1,\lambda_2,0,\ldots,0)$, where $\lambda_1 \in \mathbb{R}\setminus \{0\}$.

By Theorem~\ref{T2}, we have equality \eqref{resonance}, which in this case reads $\lambda_1+\lambda_2 = r\lambda_j$, where $j=1,2$, or $\lambda_1+\lambda_2 = 0$.
Each of these equalities defines the spectrum of $\V$ up to a common factor $\sigma$, i.e., it uniquely defines the spectrum of the corresponding direction field. In both cases $\lambda_{1,2} \in \mathbb{R}\setminus \{0\}$, hence in a neighborhood of $0$ the set of singular points of $\V$ is the center manifold $W^c$, $\codim W^c=2$.

\begin{theorem}
\label{T3-1}
Assume $W^c \subset \Gamma$, then in a neighborhood of $0$ the following statements hold.

1) There exists a smooth regular function $g:\Gamma\rightarrow \mathbb{R}$ such that $W^c=\{g=0\}$ 
and $\V{\bigl|}_{\Gamma} = g \V_1{\bigl|}_{\Gamma}$, where $\V_1{\bigl|}_{\Gamma}$ is a smooth non-vanishing field on $\Gamma$.

2) The spectrum of the linearization of $\V$ at any singular point is $\sigma (1, r-1, 0, \ldots,0)$,
where $\sigma$ is a smooth non-vanishing function on $W^c$.

3) The field $\V$ is smoothly orbitally equivalent to one of the following normal forms:
\begin{quote}
{\it

\eqref{26} with $\lambda(z) = r-1$ \ if \ $r>1$ and $r-1, \,(r-1)^{-1} \notin \mathbb N$ \ or \ $r<1$ and $r \notin \mathbb Q$,

\eqref{27} with $\lambda(z) = n$ \ if \ $r-1$ or $(r-1)^{-1}$ is equal to $n \in \mathbb N$,

\eqref{33} \ if \ $r-1 = -n/m$, where $n,m \in \mathbb N$, and $\Psi_i(0,0) \neq 0$ for at least one index $i=1,\ldots,l$ in the preliminary form \eqref{29}.

}
\end{quote}
\end{theorem}

{\bf Proof}.
For the first statement choose local coordinates $(\xi,\eta,\zeta)$ such that the invariant hypersurface $\Gamma$ is the hyperplane $\{\xi=0\}$ and the center manifold $W^c$ is the subspace $\{\xi=\eta=0\}$. Then the field $\V$ has the form
\begin{equation*}
\dot \xi = \xi v, \quad \dot \eta = \xi w_1 + \eta w_2,  \quad \dot \zeta_i = \alpha_i \xi v + \beta_i (\xi w_1 + \eta w_2), \quad  i=1,\ldots, l,
\end{equation*}
where $v, w_{1,2}$ and $\alpha_i, \beta_i$ are smooth functions of $\xi, \eta, \zeta$, and $\lambda_{1}=v(0)$, $\lambda_{2}=w(0)$. Substituting $\xi=0$, the field $\V {\bigl|}_{\Gamma}$ has the form $\dot \eta = \eta w_2$, $\dot \zeta_i = \eta \beta_i w_2$, $i=1,\ldots,l$. 
Setting the function $g=\eta w_2$, we get 
$\Gamma \cap \{g=0\} = W^c$ and $\V{\bigl|}_{\Gamma} = g \V_1{\bigl|}_{\Gamma}$, where the field $\V_1{\bigl|}_{\Gamma}$ has the form
$\dot \eta = 1$, $\dot \zeta_i = \beta_i$, $i=1,\ldots,l$.

As for the second statement, according to previous reasonings, at any singular point in a neighborhood of $0$ we have equality $\lambda_1+\lambda_2 = r\lambda_j$, where $j=1,2$, or $\lambda_1+\lambda_2 = 0$. From the hypothesis $W^c \subset \Gamma$ it follows equality \eqref{resonance} with $\lambda_j=0$ is impossible. Indeed, by Theorem~\ref{T2} the spectrum of the linearization of the restriction $\V {\bigl|}_{\Gamma}$ at $0$ is $(\lambda_1,\lambda_2,0,\ldots,0)$, where the number of zero eigenvalues is less by $1$ than in the spectrum of $\V$, i.e., is equal to $l-1$. On the other hand, the inclusion $W^c \subset \Gamma$ implies that the spectrum of the linearization 
of the restriction $\V {\bigl|}_{\Gamma}$ contains $l$ zeros.
Hence we have equality $\lambda_1+\lambda_2 = r\lambda_j$, with $j=1$ or $2$. Without loss of generality one can put $j=1$, 
then $\lambda_2 \equiv (r-1)\lambda_1$. Since the last equality holds identically at all points in $W^c$, the spectrum is $\sigma (1, r-1, 0, \ldots,0)$ with a smooth non-vanishing function $\sigma$.
The third statement follows from Theorems \ref{T2-1}\,--\,\ref{T2-4} and Remark~\ref{R1}.
\ebox

\medskip

Each of the applications in this section will cast in the following situation.

Consider the Euler--Lagrange equation 
\begin{equation}
\label{37}
\frac{d}{dt} L_p - L_x = 0,  \quad p = \frac{dx}{dt},
\end{equation}
with Lagrangian $L(t,x,p)$, where $t,x \in \mathbb R$. In the $(t,x,p)$-space equation \eqref{37} generates the direction field $\chi$ corresponding to the vector field $\W$ given by
\begin{equation}
\label{38}
\dot t = L_{pp}, \quad 
\dot x = p L_{pp}, \quad
\dot p = L_x - L_{tp} - pL_{xp},
\end{equation}
where the dot over a symbol means differentiation with respect to an independent variable playing the role of time.

\begin{lemma}
\label{L3}
At all points of the $(t,x,p)$-space where $L$ is smooth the identity $D_{\W} \equiv 0$ holds.
Consequently, at all singular points of the vector field $\W$ where $L$ is smooth, the spectrum of the linearization of $\W$ has resonance $\lambda_1+\lambda_2=0$.
The same statements are valid for the corresponding direction field $\chi$.
\end{lemma}

{\bf Proof}.
The identity $D_{\W} \equiv 0$ is due to simple calculation. The field \eqref{38} belongs to the class of vector fields of type \eqref{16}, where the generators of the ideal $I$ are $v = L_{pp}$ and $w = L_x - L_{tp} - pL_{xp}$. Hence the spectrum of the linearization of $\W$ at any singular point is $(\lambda_1, \lambda_2, 0)$. The equality $\lambda_1+\lambda_2=0$ for the field $\W$ follows from the equality $D_{\W} \equiv 0$.
The same equality for the fields $\varphi \W$ follows from the identity $D_{\varphi \W} \equiv \varphi D_{\W} + L_{\W}\varphi$.
\ebox

In the applications below we deal with the case when the Lagrangian is smooth at all points of the $(t,x,p)$-space
except for the the points of some regular surface $\Gamma=\{f=0\}$ and the components of the field $\W$ given by 
formula \eqref{38} are fractions with common denominator $f^r$, $r>0$.
Thus the field $\W$ is connected with some smooth field $\V$ by the formula \eqref{defW}. 
From the identity $D_{\W} \equiv 0$ (Lemma \ref{L3})
it follows that  conditions \eqref{conda} and \eqref{condb} will be always satisfied, hence Theorems~\ref{T1}\,--\,\ref{T3} and \ref{T3-1} are valid.

\subsection{Pseudo-Riemannian metric}

Consider a surface $S$ with a  system of coordinates $(t,x)$ and a pseudo-Riemannian metric
\begin{equation}
\label{39}
Q(dt, dx) = a(t,x) \, dx^2 + 2b(t,x) \, dx dt + c(t,x) \, dt^2
\end{equation}
with smooth coefficients $a,b,c$. The quadratic form $Q$ is positive definite on an open domain $\mathcal E \subset S$ (which is called {\it elliptic}),
indefinite on some other open domain $\mathcal H \subset S$ (which is called {\it hyperbolic}), and degenerate on the curve $A=\{\Delta=0\}$,
where $\Delta = b^2-ac$ is the discriminant of the form $Q$. The curve $A$ separates the domains $\mathcal E$ and $\mathcal H$, every point of $A$ is said {\it parabolic}.

{\bf Example 5.} Let $S$ be a smooth surface embedded in the $3$-dimensional Minkowski space, i.e., the $3$D 
affine space with Cartesian coordinates $(x,y,z)$ endowed with the pseudo-Euclidean metric $ds^2 = dx^2 + dy^2 - dz^2$. A pseudo-Riemannian metric is induced on $S$ by the metric $ds^2$ in the ambient space.
Denote by $C_P$  the {\it light cone} in the $3$D tangent space at the point $P=(x,y,z)$ given by the equation $dx^2 + dy^2 - dz^2 = 0$.
Then three possibilities arise: either the tangent plane to $S$ at $P$ does not intersect $C_P$ (then $P \in \mathcal E$), 
or it intersects $C_P$ along a pair of lines (then $P \in \mathcal H$), or finally it intersects $C_P$ along a unique line (then $P$ is parabolic).

For instance, if $S$ is a Euclidean sphere ($x^2 + y^2 + z^2 = r^2$), the parabolic points form two circles $z = \pm r/\sqrt{2}$, which separate $S$ into two elliptic domains ($\mathcal E: \, |z| > r/\sqrt{2}$) and one hyperbolic domain ($\mathcal H: \, |z| < r/\sqrt{2}$). Geodesics on Euclidean spheres and ellipsoids in 3D Minkowski space are well-studied, see e.g. \cite{GKhT}, \cite{KhT}.

\medskip

Consider geodesics generated by the pseudo-Riemannian metric \eqref{39} in a neighbourhood of a parabolic point.
Their 1-graphs are extremals of equation \eqref{37} with $L = \sqrt{F}$, where $F = a(t,x) p^2 + 2b(t,x) p + c(t,x)$.
Then the vector field $\W$ given by formula \eqref{38} reads
\begin{equation}
\label{40}
\dot t = -\Delta F^{-\frac{3}{2}}, \quad 
\dot x = -p \Delta F^{-\frac{3}{2}}, \quad
\dot p = -M F^{-\frac{3}{2}}/2,
\end{equation}
where $M = \sum\limits_{i=0}^{3} \mu_i(t,x) p^i$ is a cubic polynomial in $p$ with coefficients
\begin{multline*}
\mu_3 = a (a_t - 2b_x) + b a_x, \quad \mu_2 = b (3a_t - 2b_x) + c a_x - 2a c_x, \\
\mu_1 = b (2b_t - 3c_x) + 2c a_t - a c_t, \quad \mu_0 = c (2b_t - c_x) - bc_t.
\end{multline*}
Multiplying  $\W$  by $-F^{\frac{3}{2}}$, we obtain the field $\V$
\begin{equation}
\label{42}
\dot t = \Delta, \quad 
\dot x = p \Delta, \quad
\dot p = M/2.
\end{equation}

For any point $q_*=(t_*,x_*) \in \mathcal E \cup \mathcal H$ and any $p \in {\mathbb R}P$ 
there exists a unique geodesic passing through $q_*$ with given tangential direction $p$. 
However if $q_*$ is parabolic, this is not the case. Indeed, for any tangential direction $p \in {\mathbb R}P$ 
such that $M(q_*,p) \neq 0$ there exists a unique trajectory of  \eqref{42} passing through the point $(q_*,p)$, a vertical line, which projects onto the single point $q_*$ in the $(t,x)$-plane.
Thus, geodesics outgoing from $q_*$ must have tangential directions $p$ such that $M(q_*,p) = 0$, i.e., their  1-graphs pass through a singular point $(q_*,p)$ of the field $\V$.

Let $q_*\in A$ and consider the equation $M(q_*,p) = 0$ with respect to $p$. We shall assume  that in a neighborhood of $q_*$ the curve $A$ 
is regular and $a(q_*) \neq 0$. Then  the quadratic polynomial $F(q_*,p) = ap^2 + 2bp + c$  has a unique root $p_0(q_*) =-\frac{b(q_*)}{a(q_*)}$, that is, the isotropic direction.\footnote{The light cone at a parabolic point consists of a unique isotropic line.}
A simple substitution shows that $p_0$ is a root of the cubic polynomial $M(q_*,p)$.
Assume that the isotropic direction $p_0$ is not tangent to the curve $A$ at $q_*$, i.e, $(a\Delta_t - b\Delta_x)|_{q_*} \neq 0$.

Under the assumptions above, the cubic polynomial $M(q_*,p)$ has one or three real prime roots.\footnote{
If $S$ is a surface embedded in 3D Minkowski space, these cases correspond to positive or negative Gaussian curvature of $S$ 
in the Euclidean metric $dx^2 + dy^2 + dz^2$.
}
Define  $W^c_0 = \{q\in A, \, p=p_0(q)\}$ and  $W^c_{\pm} = \{q\in A, \, p=p_{\pm}(q)\}$ where $p_{\pm}(q_*)$ are the
non-isotropic roots of $M(q_*,p)=0$, if they exist. The union of the three curves $W^c_0$, $W^c_{\pm}$  is the set of singular points of $\V$ and coincides with its center manifold $W^c $. The function $F$ vanishes on $W^c_0$ while $F \neq 0$ on $W^c_{\pm}$. Thus the fields  \eqref{40} and \eqref{42} are connected by relation \eqref{defW}, where $f=F$ and $r=\frac{3}{2}$. 
Since the field $\W$ is obtained from an Euler--Lagrange equation,  conditions \eqref{conda} and \eqref{condb} follow from
the identity $D_{\W} \equiv 0$, which is valid for all points except for the hypersurface $\Gamma=\{F=0\}$.
From Theorem \ref{T1} it follows that $\Gamma$ is an invariant hypersurface of $\V$. Hence the isotropic curves are geodesic lines (of zero length) in the pseudo-Riemannian metric \eqref{39}. By construction $W_0^c\subset\Gamma$.

Let $(q,p_0)\in W^c_0$. Clearly, the spectrum of the linearization of $\V$ at $(q,p_0)$ contains the eigenvalue $\lambda_1 = \Delta_t + p_0\Delta_x \neq 0$.
By Theorem \ref{T3-1}, in a neighborhood of $(q,p_0)$ there exists a function $\sigma: W^c_0 \rightarrow \mathbb{R}$ such that the spectrum of the linearization of $\V$ at all points sufficiently close to $(q,p_0)$ is $\sigma (2,1,0)$. Computing, 
we easily get $\sigma = \Delta_t + p_0\Delta_x$.
Hence the germ of $\V$ at $(q,p_0)$ is smoothly orbitally equivalent to
\begin{equation}
\dot \xi = 2\xi + \varphi(\zeta) \eta^2, \quad \dot \eta = \eta, \quad \dot \zeta = 0.
\label{43}
\end{equation}
The normal form \eqref{43} can be further simplified.

\begin{theorem}
\label{T3-2}
The germs of the vector field $\V$ given by formula \eqref{42} at the singular points $(q,p_0) \in W^c_0$ and $(q,p_{\pm}) \in W^c_{\pm}$ are  smoothly orbitally equivalent to
\begin{eqnarray}
\label{44}
&& \dot \xi = 2\xi, \quad \dot \eta = \eta, \quad \phantom{...} \dot \zeta = 0, \\
\label{45}
&& \dot \xi = \xi, \quad \phantom{..} \dot \eta = -\eta, \quad \dot \zeta = \xi \eta,
\end{eqnarray}
respectively.
\end{theorem}

{\bf Proof}.
To establish normal form \eqref{44} it is sufficient to prove that the coefficient $\varphi(\zeta)$ in the normal form \eqref{43} is identically equal to zero. Indeed, the field \eqref{43} has the invariant foliation $\{(\xi,\eta,\zeta) : \zeta=c\}$, and the restriction to each leaf is a node with exponent\footnote{
The exponent of a node (or saddle) is defined to be the ratio of the eigenvalue of largest modulus of the linearization
field to the smallest one.
}
equal to $2$. The eigenvalue of largest modulus corresponds to the eigenvector $\frac{\partial}{\partial \xi}$ and the eigenvalue of smallest modulus corresponds to $\frac{\partial}{\partial \eta}$.

Given an arbitrary point $(q_*,p_0)\in W^c_0$, consider the restriction of the field \eqref{43} to the invariant leaf $\{(\xi,\eta,\zeta):\zeta=\zeta_*\}$ passing through $(q_*,p_0)$. Integrating the corresponding differential equation $d\xi/d\eta = 2\xi/\eta + \varphi(\zeta_*)\eta$, we get the single integral curve $\eta = 0$ and the family of integral curves
\begin{equation}
\xi = c\eta^2 + \varphi(\zeta_*)\eta^2 \ln |\eta|, \quad c={\rm const},
\label{family}
\end{equation}
with common tangential direction $\frac{\partial}{\partial \eta}$ at $0$. 

In the case $\varphi(\zeta_*) = 0$ all curves of the family \eqref{family} are parabolas,
in the case $\varphi(\zeta_*) \neq 0$ they are $C^1$-smooth, but not $C^2$-smooth at $0$.
On the other hand, the previous reasoning shows that the germ of $\V$ at $(q_*,p_0)$ has at least one $C^{\infty}$-smooth integral curve: the vertical line (parallel to the $p$-axis). Simple calculation shows that the direction $\frac{\partial}{\partial p}$ in the initial coordinates $(t,x,p)$ corresponds to the direction $\frac{\partial}{\partial \eta}$ in the normal coordinates $(\xi,\eta,\zeta)$. Hence family \eqref{family} contains at least one $C^{\infty}$-smooth integral curve. This implies that $\varphi(\zeta_*) = 0$.

The second statement of the theorem (the normal form \eqref{45}) follows from Lemma \ref{L3} and Theorem \ref{T2-4}; validity of the condition $\Psi_1(0,0)\neq 0$ can be proved by direct calculation (see Theorem~2 in \cite{Rem-Pseudo}).
\ebox

\subsection{Metric of  Klein type}

A natural generalization of the Klein metric on the $(t,x)$-plane is
\begin{equation}
ds^2 = \frac{\alpha \, dx^2 + 2\beta \, dx dt + \gamma \, dt^2}{t^{2n}}, \quad n \in \mathbb N,
\label{46}
\end{equation}
where the numerator is a positive definite quadratic form\footnote{
The case when the numerator is an indefinite (and non-degenerate) quadratic form was also studied \cite{Rem-Klein},
but for our present purposes it is sufficient to consider  the positive definite case.
}
with coefficients $\alpha, \beta, \gamma$ smoothly  depending on $t,x$. We study locally the geodesics of metric \eqref{46} passing through a singular point, i.e., a point of the axis $\{t=0\}$. It is not hard to prove that in appropriate local coordinates on the $(t,x)$-plane the germ of metric \eqref{46} simplifies to the form
\begin{equation}
ds^2 = \frac{\alpha \, dx^2 + \gamma \, dt^2}{t^{2n}}, \quad n \in \mathbb N,
\label{47}
\end{equation}
with smooth positive coefficients $\alpha(t,x)$ and $\gamma(t,x)$.

The geodesics of metric \eqref{47} are extremals of the Euler--Lagrange equation \eqref{37} with $L = \sqrt{F}/t^n$, where $F = \alpha p^2 + \gamma > 0$ and $p=dx/dt$. The corresponding vector field $\W$ reads
\begin{equation}
\dot t =  \alpha \gamma t^{-n} F^{-\frac{3}{2}}, \quad 
\dot y = \alpha \gamma p t^{-n} F^{-\frac{3}{2}}, \quad 
\dot p = - \frac{1}{2} t^{-n-1} M F^{-\frac{3}{2}}, 
\label{48}
\end{equation}
where $M = \sum\limits_{i=0}^{3} \mu_i(t,x) p^i$ is a cubic polynomial of $p$ with coefficients
\begin{equation*}
\mu_3 = \alpha (t\alpha_t - 2n \alpha),  \quad \mu_2 = t (\alpha_x \gamma - 2\alpha \gamma_x), \quad
\mu_1 = t (2\alpha_t \gamma - \alpha \gamma_t) - 2n \alpha \gamma, \quad \mu_0 = -t\gamma \gamma_x.
\end{equation*}
Multiplying $\W$ by $f^{n+1}$, where $f =tg^{\frac{1}{n+1}}$ and $g=F^{3/2}/(\alpha \gamma) > 0$, we obtain the field $\V$
\begin{equation}
\label{50}
\dot t = t, \quad \dot x = pt, \quad \dot p = -M/(2\alpha \gamma).
\end{equation}

Fields \eqref{48} and \eqref{50} are connected by relation \eqref{defW}, where $f=tg^{\frac{1}{n+1}}$ and $r=n+1$.
Conditions \eqref{conda} and \eqref{condb} are satisfied, $\W$ being obtained from an Euler--Lagrange equation.
Theorem \ref{T1} implies that $\Gamma = \{f=0\}= \{t=0\}$ is an invariant plane for $\V$.
The restriction of the field $\V$ to $\Gamma$ is parallel to the $p$-axis. Hence geodesics outgoing from a point $q_*=(0,x_*)$ must have tangential directions corresponding to $p$ such that $M(q_*,p) = 0$, i.e., their 1-graphs pass through singular points of $\V$.

Given a point $q_*=(0,x_*)$, consider the equation $M(q_*,p) = 0$ with respect to $p$. Since $M(q_*,p) = -2n\alpha p (\alpha p^2 + \gamma)$, the cubic polynomial $M(q_*,p)$ has the only real root $p=0$. The spectrum of the linearization of $\V$ at  $(q_*,0)$ is $ (n,1,0)$, and the $x$-axis is the center manifold ($W^c$). Clearly,  $W^c \subset \Gamma$ and from Theorem \ref{T3-1} we get the following result.

\begin{theorem}
\label{T3-3}
The germ of the vector field \eqref{50} at the singular point $(q_*,0)$ is smoothly orbitally equivalent to 
\begin{equation}
\label{51}
\dot \xi = n\xi + \varphi(\zeta)\eta^n, \quad  \dot \eta = \eta, \quad \dot \zeta = 0.
\end{equation}
\end{theorem}

Unlike the case of geodesics in pseudo-Riemannian metrics, here the coefficient $\varphi(\zeta)$ is not necessarily zero.
For instance, in the case $n=1$ the condition $\varphi(\zeta_*) = 0$ is equivalent to $\gamma_x(0,x_*) = 0$, where $\zeta_*$ corresponds to the point $q_*=(0,x_*)$. Clearly, if $\varphi(0) \neq 0$ then $\varphi(\zeta)$ simplifies to~$1$, if $\varphi(\zeta)$ has a finite order $s$ at the origin then $\varphi(\zeta)$ simplifies to $\zeta^{s}$.

{\bf Example 6.}
Consider the Klein metric, given by formula \eqref{47} with $\alpha \equiv \gamma \equiv 1$ and $n=1$. In this case vector field $\V$ given by \eqref{50} has the normal form \eqref{51} with $\varphi(\zeta) \equiv 0$, since $\gamma_x(0,x) \equiv 0$. Hence the restriction of the field $\V$ on each invariant leaf (given by the formula $\zeta = c$ in the normal coordinates) is a bicritical node. Thus the integral curves of $\V$ are $C^{\infty}$-smooth, and for each singular point $q_*=(0,x_*)$ there exists a family  of geodesics with common tangential directions $p=0$ and various 2-jets.
Indeed, geodesics of the Klein metric passing through the point $q_* \in A$ (here $A$ is the absolute) are the circles $(x-x_*)^2 + t^2 = R^2$ and the straight line $x=x_*$.

\subsection{Almost-Riemannian metric}

Let $\V_1, \V_2$ be smooth vector fields on the $(x,y)$-plane. Assume them to be collinear at the points of a regular curve $A$. The metric $ds^2$ having $(\V_1,\V_2)$ as orthonormal frame is well-defined, smooth and positive definite on the whole plane except for the points of  $A$. Our aim is to study  geodesics of the metric $ds^2$ passing through a point of the curve~$A$.

According to \cite{ABS}, for a generic pair $\V_1, \V_2$ in a neighborhood of almost all points of $A$ there exist local coordinates such that
\begin{equation}
\V_1 = \frac{\partial}{\partial x}, \quad \V_2 = 2x v^{-1}(x,y) \frac{\partial}{\partial y},
\label{new}
\end{equation}
where $v(x,y)$ is a smooth non-vanishing function.
Points in $A$ at which such a coordinate system does not exist form a discrete subset of $A$ and will not
be considered in the following.
The required metric for  fields \eqref{new} is
\begin{equation*}
ds^2 = dx^2 + \frac{v^2}{4x^2} \, dy^2 = \frac{v^2 \, dy^2 + d(x^2)^2}{4x^2}.
\end{equation*}
Substituting $t=x^2$ and multiplying by the unessential constant factor $4$, we get
\begin{equation}
ds^2 = \frac{{\widetilde v}^2 \, dy^2 + dt^2}{t},  \quad {\widetilde v}(t,y)=v_1(t,y) + \sqrt{t}\,v_2(t,y),
\label{53}
\end{equation}
where $v_{1,2}$ are smooth functions defined by decomposition of the function $v(x,y)$ into even and odd parts:
$v(x,y)=v_1(x^2,y) + x v_2(x^2,y)$. Geodesics of the metric \eqref{53} are extremals of Euler--Lagrange equation with Lagrangian $L = \sqrt{F/t}$, where $F = {\widetilde v}^2 p^2 + 1$ and $p=dy/dt$. The corresponding vector field $\W$ in the $(t,y,p)$-space reads
\begin{equation}
\dot t =  {\widetilde v}^2 t^{-\frac{1}{2}} F^{-\frac{3}{2}}, \quad 
\dot y =  {\widetilde v}^2 p t^{-\frac{1}{2}} F^{-\frac{3}{2}}, \quad 
\dot p =  \frac{{\widetilde v}}{2} t^{-\frac{3}{2}} F^{-\frac{3}{2}} {\widetilde M},
\label{55}
\end{equation}
where ${\widetilde M} = \sum\limits_{i=0}^{3}\,{\widetilde \mu}_i(t,y) p^i$ is a cubic polynomial in $p$ with coefficients
\begin{equation}
{\widetilde \mu}_3 = {\widetilde v}^3 - 2t{\widetilde v}^2 {\widetilde v}_t, \quad  {\widetilde \mu}_2 = -2t {\widetilde v}_y, \quad {\widetilde \mu}_1 = {\widetilde v} - 4t {\widetilde v}_t, \quad {\widetilde \mu}_0 = 0.
\label{56}
\end{equation}

Multiplying $\W$ by $f^r$, where $r=\frac{3}{2}$, $f=tg$ and $g=(2/{\widetilde v})^{\frac{2}{3}}F \neq 0$, we obtain the field $\V$
\begin{equation}
\label{57}
\dot t = 2{\widetilde v} t, \quad \dot y = 2{\widetilde v} tp, \quad \dot p = {\widetilde M}.
\end{equation}
Fields \eqref{55} and \eqref{57}  are connected by relation \eqref{defW}, where the function $f=tg$ is regular and $r=\frac{3}{2}$. Nevertheless, in general we cannot apply Theorems~\ref{T1}, \ref{T2} and \ref{T3-1}, since the field $\V$ is not even $C^1$-smooth. Indeed, the components of the field $\V$ depend on the function ${\widetilde v}(t,y)$ and its first-order derivatives, which are smooth only if $v(x,y)$ is an even function of $x$ (see formula \eqref{53}).

{\bf Example 7.}
Consider the Clairaut--Liouville metric. This is an example in which the vector field \eqref{57} turns out to be smooth, the function $v(x,y)$ being even in $x$. For instance, in \cite{BCT} the authors deal with the metric
\begin{equation*}
ds^2 = dx^2 + \frac{g(x^2,y)}{x^2} \, dy^2 = \frac{x^2 \, dx^2 + g(x^2,y)\, dy^2}{x^2},
\end{equation*}
where $g$ is a positive smooth function ($x$ and $y$ are standard angular coordinates on the sphere, the curve $A=\{x=0\}$ is the equator).\footnote{
In the case $g(x^2,y) \equiv 1$ this formula gives the well-known Grushin metric.
}
After the change of variables $t=x^2$ we get the metric \eqref{53} with ${\widetilde v} = 2\sqrt{g(t,y)}$, which leads to the smooth field \eqref{57}.

\medskip

To overcome the problem, we make the change of variable $x^2=t$ in \eqref{57}. This yields to
\begin{equation}
\label{58}
\dot x = xv, \quad \dot y = 2 x^2 vp , \quad \dot p = M,
\end{equation}
where $M = \sum\limits_{i=0}^{3} \mu_i(x,y) \,p^i$ and $\mu_i(x,y)={\widetilde \mu}_i(x^2,y)$. The coefficients ${\widetilde \mu}_i$ are polynomials of the function ${\widetilde v}(t,y)$ and its first-order derivatives (see formulas \eqref{56}). Note also that ${\widetilde v}_t$ appear in \eqref{56} with the factor $t$, whence after the substitution $x^2=t$ the expression $t{\widetilde v}_t$ becomes a smooth function of $x,y$.

The first two components of the field \eqref{58} vanish at $x=0$. Given a point $q=(0,y)$ consider the cubic equation $M(q,p) = 0$ with respect to $p$.
It reads $v(q) p ((v(q) p)^2 + 1) = 0$. This equation has a unique real root $p_0=0$.
Recalling that $p=dy/dt$, the root $p_0=0$ defines the unique admissible direction for geodesics passing through the point $(0,y)$ on the $(t,y)$-plane.
The corresponding  direction on the $(x,y)$-plane is given by the relation $dy/dx=2xp$ which is also equal to zero.
The spectrum of the linearization of the germ \eqref{58} at $(q,p_0)$ is $(\lambda_1,\lambda_2,0)$, where $\lambda_1=v(q)$ and
$\lambda_2=M_p(q,p_0) = v(q)$.

\begin{theorem}
\label{T3-4}
The germ of the vector field \eqref{58} at the singular point $(q,p_0)$ is smoothly orbitally equivalent to
\begin{equation}
\label{59}
\dot \xi = \xi, \quad  \dot \eta = \eta, \quad \dot \zeta = 0.
\end{equation}
\end{theorem}

{\bf Proof}.
By Theorem~\ref{T2-3}, the germ of the vector field \eqref{58} at $(q,p_0)$ is smoothly orbitally equivalent to normal form \eqref{27} with $\lambda(\zeta) \equiv 1$ and $l=1$. To establish normal form \eqref{59} it is sufficient to prove that the coefficient $\varphi(\zeta)$ in \eqref{27} is identically equal to zero.

Let $\Lambda$ be the linearization of the vector field \eqref{58} at the singular point $(q,p_0)$. Consider the matrix $\Lambda-\lambda I$, where $\lambda=v(q)$ is the double eigenvalue of $\Lambda$. Clearly, the value $\rg (\Lambda-\lambda I)$ equals either 1 or 2 and it is an invariant of the field.
Hence $\varphi(0)=0$ if $\rg (\Lambda-\lambda I)=1$ and $\varphi(0)\neq 0$ if $\rg (\Lambda-\lambda I)=2$.
On the other hand, a simple calculation shows that $\rg (\Lambda-\lambda I)=1$ if and only if $M_x(q,p_0) = 0$.
Recalling that $M(q,p) = v(q) p ((v(q) p)^2 + 1)$ and $p_0=0$ we get $M_x(q,p_0) = 0$.
This completes the proof.
\ebox

From the normal form \eqref{59} it follows that vector field \eqref{58} has an invariant foliation (given by $\zeta = {\rm const}$ in the normal coordinates) such that each leaf intersects the center manifold $W^c$ at a unique point. Hence, geodesics passing through the point $q_*=(0,y_*)$ on the $(x,y)$-plane are projections of integral curves lying in the corresponding leaf. The restriction of vector field \eqref{58} to the leaf is a bicritical node, hence there is a one-parameter family of integral curves passing through the point $q_*=(0,y_*)$. This gives a family of smooth geodesics passing through the point $q_*$ with common tangential direction which coincides with $\V_1$ and $\V_2$ at the point $q_*$.
Moreover, the geodesics have the same 2-jet and different 3-jets at $q_*$.

{\bf Example 8.}
Geodesics in the Grushin metric (which corresponds to the vector fields \eqref{new} with $v(x,y)\equiv 2$) have the form
$y(x) = y_* + c^{-2} \arcsin (cx) - c^{-1}x \sqrt{1-c^2x^2}$, where $c$ is an arbitrary constant.


\end{document}